\documentclass{amsart}

\usepackage[utf8]{inputenc}
\usepackage{amsmath}
\usepackage{latexsym}
\usepackage{amsfonts,amssymb}
\usepackage{color} \usepackage{mathrsfs} \let\mathcal\mathscr
\usepackage[all,ps,cmtip]{xy}
\usepackage{lscape}\usepackage{amsthm}
\usepackage{fancyhdr}
\usepackage[T1]{fontenc}
\usepackage{graphicx}
\usepackage{tikz}
\usepackage{eurosym}
\usepackage{skull}
\usepackage{mathenv}
\usepackage{pstricks}

\title{Gradings of Lie algebras, magical spin geometries  and matrix factorizations}
\author{Roland Abuaf, Laurent Manivel}
\date{January 2019}
\DeclareGraphicsExtensions{eps}
\DeclareGraphicsExtensions{pdf}

\newtheorem{theo}{Theorem}[subsection]
\newtheorem{theo*}{Theorem}

\newtheorem{prop}[theo]{Proposition}
\newtheorem{prop*}{Proposition}

\newtheorem{prob*}{Problem}

\newtheorem{lem}[theo]{Lemma}

\newtheorem{theos}{Th\'eor\`eme}[section]

\newtheorem{rems}[theos]{Remarks}

\def\DB{\mathrm{D^{b}}}
\def\OO{\mathcal{O}}
\def\w{\omega}
\def\DM{\mathrm{D^{-}}}
\def\DP{\mathrm{D^{perf}}}
\def\Ri{\mathrm{R^{i}}}
\def\R0{\mathrm{R^{0}}}
\def\ot{\otimes}
\def\Hi{\mathrm{H^i}}
\def\f{\tilde{f}}
\def\HH{\mathrm{HH}}
\def\Hh{\mathcal{H}om}
\def\HHH{\mathrm{Hom}}
\def\LL{\mathrm{\textbf{L}}}
\def\Li{\mathrm{L^{i}}}
\def\RR{\mathrm{\textbf{R}}}
\def\OO{\mathbb{O}}\def\AA{\mathbb{A}}\def\BB{\mathbb{B}}
\def\RR{\mathbb{R}}\def\CC{\mathbb{C}}\def\HH{\mathbb{H}}
\def\cE{\mathcal{E}}
\def\PP{\mathbb{P}}
\def\D{\mathcal{D}}
\def\d{\delta}
\def\tY{\tilde{Y}}
\def\tU{\tilde{U}}
\def\tV{\tilde{V}}
\def\tF{\tilde{F}}
\def\P{\mathcal{P}}
\def\M{\mathcal{M}}
\def\C{\mathcal{C}}
\def\F{\mathcal{F}}
\def\E{\mathcal{E}}
\def\GG{\mathrm{G}}
\def\G{\mathcal{G}}
\def\H{\mathcal{H}}
\def\Q{\mathcal{Q}}
\def\T{\mathcal{T}}
\def\A{\mathcal{A}}
\def\cL{\mathcal{L}}\def\cO{\mathcal{O}}
\def\cE{\mathcal{E}}\def\cV{\mathcal{V}}
\def\Sp{\mathrm{Sp}}
\def\B{\mathcal{B}}
\def\hX{\hat{X}}
\def\tX{\tilde{X}}
\def\X{\mathcal{X}}
\def\J{\mathrm{J}}
\def\w{\omega}
\def\vT{\mathrm{vol}^*(\T)}
\def\spin{\mathrm{Spin}}\def\sl{\mathrm{SL}}
\def\V{\mathrm{V}_7}
\def\U{\mathrm{U}_6}
\def\SS{\mathrm{S}}
\def\ZZ{\mathbb{Z}}
\def\fg{\mathfrak{g}}\def\ft{\mathfrak{t}}\def\fh{\mathfrak{h}}
\def\fs{\mathfrak{s}}\def\fe{\mathfrak{e}}\def\fsp{\mathfrak{sp}}
\def\fsl{\mathfrak{sl}}\def\fso{\mathfrak{so}}\def\ff{\mathfrak{f}}
\def\fspin{\mathfrak{spin}}\def\fgl{\mathfrak{gl}}
\def\lra{\longrightarrow}\def\ra{\rightarrow}

\DeclareMathOperator{\Z}{\mathbb{Z}}
\DeclareMathOperator{\sing}{Sing}
\DeclareMathOperator{\ad}{ad}
\DeclareMathOperator{\codim}{Codim}
\DeclareMathOperator{\Hom}{Hom}
\DeclareMathOperator{\sym}{Sym}
\DeclareMathOperator{\Pic}{Pic}
\DeclareMathOperator{\rg}{rg}
\DeclareMathOperator{\Gr}{Gr}
\DeclareMathOperator{\ch}{ch}
\DeclareMathOperator{\Pff}{Pff}

\renewcommand{\refname}{References}

\newcommand{\eq}[1][r]
{\ar@<-3pt>@{-}[#1]
\ar@<-1pt>@{}[#1]|<{}="gauche"
\ar@<+0pt>@{}[#1]|-{}="milieu"
\ar@<+1pt>@{}[#1]|>{}="droite"
\ar@/^2pt/@{-}"gauche";"milieu"
\ar@/_2pt/@{-}"milieu";"droite"}

\newcommand{\incl}[1][r]
  {\ar@<-0.2pc>@{^(-}[#1] \ar@<+0.2pc>@{-}[#1]}

\newcommand{\leftexp}[2]{{\vphantom{#2}}^{#1}{#2}}

\begin{document}

\maketitle

\begin{abstract}
We describe a remarkable rank $14$ matrix factorization of the octic $\spin_{14}$-invariant polynomial on either of its half-spin 
representations. We observe that this representation can be, in a suitable sense, identified with a tensor product of two octonion 
algebras. Moreover the matrix factorisation can be deduced from a particular $\ZZ$-grading of $\fe_8$. Intriguingly, the whole story 
can in fact be extended to the whole Freudenthal-Tits magic square and yields matrix factorizations on other spin representations,
as well as for the degree seven invariant on the space of three-forms in several variables. As an application of our results on 
$\spin_{14}$, we construct
a special rank seven vector bundle on a double-octic threefold, that we conjecture to be spherical.
%Finally, we use the same circle of ideas to exhibit matrix factorizations of the semi-invariants 
%of degree seven and sixteen on the spaces of three-forms in seven and eight variables.
\end{abstract}

%\vspace{\stretch{1}}

%\newpage

%\tableofcontents

\section{Introduction}

Recall that a {\it matrix factorization} of a polynomial $W$ is  a pair $(P,Q)$ 
of square matrices of the same size, say $N$, with polynomial entries,  such that 
$$PQ=QP=W.Id_N.$$
Matrix factorizations have attracted a lot of attention since their introduction by 
Eisenbud \cite{eisenbud} in connection with Cohen-Macaulay modules over hypersurfaces. 
Important examples of matrix factorizations, when $W$ is a quadratic form, are provided 
by Clifford modules \cite{beh, bertin}. They can be obtained as follows. Suppose our 
base field is the field of complex numbers, and consider the simple Lie algebras $\fso_n$,
$n\ge 5$, with their spin representations. When $n$ is even, there are two half-spin 
representations $\Delta_+$ and $\Delta_-$, of the same dimension $N=2^{\frac{n}{2}-1}$. 
Their direct sum can be defined as a module
over the Clifford algebra of the natural representation $V_n$ of $\fso_n$, with its invariant 
quadratic form $q$. The Clifford multiplication yields equivariant morphisms 
$$V_n\otimes\Delta_+\lra \Delta_-\qquad\mathrm{and}\qquad V_n\otimes\Delta_-\lra \Delta_+.$$
So for each $v\in V_n$, we get morphisms $P(v): \Delta_+\ra \Delta_-$ and  $Q(v): \Delta_-\ra \Delta_+$, 
depending linearly on $v$, and the fact that the total spin representation is a Clifford module 
yields the identities 
$$P(v)\circ Q(v)= q(v)Id_{\Delta_-} \qquad\mathrm{and}\qquad Q(v)\circ P(v)=q(v)Id _{\Delta_+}.$$
In other words, we get a of rank $N$ matrix factorization of the quadratic form $q$. 

Surprisingly, this is a non trivial matrix factorization of minimal size of a non degenerate
quadratic form in $n$ variables. 
This illustrates the difficulty to find explicit ones in general. One of the goals 
of this paper is precisely to describe several remarkable matrix factorizations, again related to spin 
representations. Our main result will be the description of a rank $14$ matrix factorization 
of a particular degree eight polynomial in $64$ variables, a $\spin_{14}$-invariant polynomial
on a half-spin representation $\Delta_{14}$. In the next section, we will give a direct construction 
of this invariant and prove that it admits a matrix factorization (Theorem \ref{mf14}). Our proof is 
remarkably simple, and relies on the fact that $\Delta_{14}$ contains an open orbit of the action 
of $\CC^*\times \spin_{14}$. We will observe in passing the intriguing fact that fixing a point in this 
open orbit determines a factorisation of $\Delta_{14}$ as the tensor product of two octonion algebras
(Proposition \ref{octooctonions}).

\smallskip
In the last section, we will relate those observations to gradings of Lie algebras and the 
Freudenthal magic square. The point is that $\Delta_{14}$ appears in a particular $\ZZ$-grading
of the largest exceptional algebra $\fe_8$, and that the octic invariant and its matrix 
factorization can be constructed directly from  $\fe_8$. Moreover, this particular $\ZZ$-grading
turns out to be related to the space of "points" in the Freudenthal geometry associated to $E_8$.

Astonishingly, the whole story extends to the full magic square. Recall that  this square
associates to a pair $(\AA,\BB)$ of normed algebras (either $\RR, \CC$, the algebra $\HH$ of 
quaternions or the Cayley algebra $\OO$ of octonions)
a semisimple Lie algebra $\fg(\AA,\BB)$. In particular $\fg(\OO,\OO)=\fe_8$. 
The space of "points" in the corresponding Freudenthal geometry induces a grading of $\fg(\AA,\BB)$ 
whose main component is always a spin representation,  and this yields a matrix factorization 
(Theorem \ref{magicmf}). Moreover, once one chooses a general point in that 
representation, it gets naturally identified with $\AA\otimes\BB$. 

Finally, we discuss the sporadic case of the third exterior power of a vector space of 
dimension seven, which is related with a certain $\ZZ$-grading of $\fe_7$.
%respectively. Moreover they admit semi-invariants of degree seven
%and sixteen, for which we obtain matrix factorizations of rank seven and eight, respectively 
%(Theorem \ref{hepticmf} and Theorem \ref{16mf}). 

\smallskip 
One motivation for this study of matrix factorizations has been the construction of a special rank seven vector bundle on a 
double octic threefold obtained as a double cover of $\PP^3$ branched over a linear section 
of the octic hypersurface defined by the $\spin_{14}$-invariant of $\Delta_{14}$. We conjecture 
that this bundle is spherical. Such a double cover is in fact a Calabi-Yau threefold and there is
an astonishing series of relationships, far from being completely understood yet, between 
exceptional Lie algebras and certain families of manifolds of Calabi-Yau type \cite{abuaf, im}. 
We hope to come back to this conjecture in a subsequent paper. 

\smallskip\noindent  {\bf Acknowledgments}. We warmly thank Vladimiro Benedetti for his help with LiE.

\section{Spin geometry in dimension fourteen} 

Spin geometry in dimension twelve has several very remarkable features, two of which
we would like to recall briefly. Let $\Delta_{12}$ be one of the half-spin representations
of $\spin_{12}$ (see \cite[Section 5]{LMfreud} for more details). 
\begin{enumerate}
\item The action of $\spin_{12}$ on $\PP\Delta_{12}$ has only four orbits, whose closures
are the whole space, a degree four hypersurface, its singular locus, and inside the latter,
the spinor variety $S_{12}$, which parametrizes one of the families of maximal isotropic
subspaces of a quadratic twelve dimensional vector space $V_{12}$.
\item The spinor variety $S_{12}\subset\PP\Delta_{12}$ is a variety {\it with one 
apparent double point}, which means that through a general point of $\PP\Delta_{12}$
passes a unique line which is bisecant to $S_{12}$. One can deduce
that the open orbit $\cO_0\simeq \spin_{12}/\sl_6\rtimes\ZZ_2$, where the $\ZZ_2$ factor
exchanges the two points in $S_{12}$ of the afore mentionned bisecant. 
\end{enumerate}
An interesting consequence is that a double cover of the open orbit is naturally built in 
the spin geometry, which turns out to be intimately related with the family of double
quartic fivefolds. Those varieties have attracted some interest from the early ages of
mirror symmetry, being Fano manifolds of Calabi-Yau type that can be considered as mirror
to certain rigid Calabi-Yau threefolds \cite{cdp, sch}. 

The goal of this section is to describe the similar properties that can be observed for the 
spin geometry in dimension fourteen. We will start by briefly recalling the orbit structure,
which has classially been considered by several authors \cite{sk, popov, kv}. The quartic invariant 
hypersurface in $\PP\Delta_{12}$ is in particular replaced by an octic invariant 
hypersurface in $\PP\Delta_{14}$ on which the next section will focus. Here we will highlight 
a kind of multiplicative version of the one apparent double point property, which we find 
very remarkable.

\subsection{Orbits} Recall that the half-spin representations of $\spin_{14}$ can be
defined by choosing a splitting $V_{14}=E\oplus F$, where $E$ and $F$ are maximal
isotropic subspaces. Note that the quadratic form $q$ on $V_{14}$ induces a perfect
duality between $E$ and $F$. Then $E$ acts on the exterior algebra $\wedge^* E$ by the 
wedge product, and $F$ by twice the contraction by the quadratic form. 
The resulting action of $V_{14}$ on $\wedge^* E$
upgrades to an action of its Clifford algebra. By restriction, on gets an action of $\spin_{14}$, as 
well as of its Lie algebra $\fspin_{14}\simeq\wedge^2V_{14}$. The half-spin representations are then 
given by the even and odd parts, $\wedge^+E$ and $\wedge^-E$,  of the exterior algebra 
$\wedge^*E$. We will let $\Delta_{14}=\wedge^+ E$. 
(Of course the construction works to any $\spin_{2n}$, starting from a splitting of a $2n$-dimensional
vector space $V_{2n}$ endowed with a non degenerate quadratic form. In the odd case, $V_{2n+1}$ can 
only be split as $E\oplus F\oplus L$, with $E$ and $F$ isotropic and $L$ a line. The unique
spin representation can then be identified with $\wedge^*E$. See e.g. \cite{chevalley} for more details.)

According to Sato and Kimura \cite[page 132]{sk}, the fact that a half-spin representation 
$\Delta_{14}$ of $\spin_{14}$ is prehomogeneous under the action of $\CC^*\times \spin_{14}$ 
was first observed by Shintani in 1970, and the orbit structure was obtained by Kimura and
Ozeki in 1973. The fact that there are only finitely many orbits is actually an immediate 
consequence of Kac and Vinberg's theory of $\theta$-groups \cite{kac, vinberg}. 
Indeed, the half-spin representation $\Delta_{14}$ 
is a component of the $\ZZ$-grading of $\fe_8$ defined by its first simple root. 

Up to our knowledge, the classification of the orbits of  $\spin_{14}$ was first published 
in 1977 by Popov \cite{popov}, with explicit representatives of each orbit and the types of 
their stabilizers. It also appears in the paper  by Kac and Vinberg 
in \cite{kv}, along with the the orbits of  $\spin_{13}$ on the same representation. 
More details about the geometry of the orbit closures can be found in \cite{kwE8}. 

As we already mentionned there is an octic invariant $J_8$ (unique up to scalar), and each level set 
$J_8^{-1}(c)$ is a single orbit of 
$\spin_{14}$ for $c\ne 0$. Inside the octic hypersurface $(J_8=0)$, 
there are eight non trivial orbits. Among those, the most important one is the (pointed) cone 
over the spinor variety $S_{14}$, which parametrizes the maximal isotropic subspaces of 
$V_{14}$ in the same family as $F$. The other family $S_{14}'$ of such spaces, to which $E$ belongs, 
is naturally embedded inside the projectivization of the other half-spin representation, the dual
$\Delta_{14}^\vee$. By the way, although we will not use this fact, it is interesting to note that 
the projective dual of $S_{14}'$ is precisely the octic hypersurface 
$(J_8=0)$ inside $\PP(\Delta_{14})$. 

For future use let us mention that the following points $z_0$ and $z_1$ are respectively 
outside the octic hypersurface, and inside the open orbit of this hypersurface \cite{kv}:
\begin{equation}\label{basepoint0}
z_0 =  1+e_{1237}+e_{4567}+e_{123456}, 
\end{equation}
\begin{equation}\label{basepoint1}
z_1 =  1+e_{1237}+e_{1587}+e_{2467}+e_{123456}.
\end{equation}
Here we have fixed a basis $e_1,\ldots , e_7$ of $E$, and we have used the notation $e_{ijkl}$ 
for $e_i\wedge e_j\wedge e_k\wedge e_l$, and so on. We will also denote by $f_1,\ldots , f_7$
the dual basis of $F$, defined by the condition that $q(e_i,f_j)=\delta_{ij}$. 

The stabilizer of $z_0$ was computed explicitely  in \cite[Proposition 40]{sk}, or more 
precisely its Lie algebra $\fh\simeq \fg_2\times\fg_2$. Sato and Kimura observed that 
$\fh$ stabilizes a unique pair of seven-dimensional subspaces of $V_{14}$, each copy of 
$\fg_2$ acting irreducibly on one of these spaces and trivially on the other. 
To be more specific,  the two
invariant seven-dimensional subspaces are  
$$V_7=\langle e_1,e_2,e_3,f_1,f_2,f_3,e_7-f_7\rangle \quad  \mathrm{and}\quad 
V'_7=\langle e_4,e_5,e_6,f_4,f_5,f_6,e_7+f_7\rangle.$$
Observe that $V_7$ and $V_7'$ are mutually orthogonal, and in direct sum (equivalently,
the restriction of the quadratic form to these spaces is non degenerate). Moreover 
$E_3:=E\cap V_7=\langle e_1,e_2,e_3\rangle$ and $E'_3:=E\cap V'_7=\langle e_4,e_5,e_6\rangle$
are maximal isotropic subspaces of $V_7$ and $V'_7$. In particular we may identify the 
spin representations of $\fspin(V_7)$ and $\fspin(V'_7)$ with $\wedge^*E_3$ and $\wedge^*E'_3$,
respectively. Observe that their tensor product is 
$$\wedge^*E_3\otimes\wedge^*E'_3 \simeq \wedge^*(E_3\oplus E'_3)\simeq \wedge^+E=\Delta_{14}.$$
Under this identification, 
$$z_0=1\otimes 1+e_{123}\otimes 1+1\otimes e_{456}+e_{123}\otimes e_{456}=
(1+e_{123})\otimes (1+e_{456})$$
turns out to be decomposable. Moreover, the action of $\CC^*\times \spin_7$ on the spin 
representation is prehomogeneous; in fact there is a quadratic invariant $Q$, and the 
non trivial orbits are the pointed quadric and its complement; the vector $1+e_{123}$
belongs to the latter, and its stabilizer is isomorphic to 
$G_2=Aut(\mathbb{O})$ (by \cite[Proposition 25]{sk}, the generic stabilizer has Lie algebra
$\fg_2$; that it is really isomorphic to $G_2$ follows from the triality 
principle \cite[Section 2.4]{baez}). 

Any element $g$ in the stabilier of $z_0$ which fixes $V_7$ must satisfy 
$g(1+e_{123}) = 1+e_{123}$, and similarly for $V'_7$. We can therefore 
conclude that $G_2\times G_2$ is the subgroup of this stabilizer that 
fixes $V_7$ and $V'_7$. The stabilizer itself will be bigger only if it contains a 
transformation that swaps those two spaces, and this is indeed what happens. 

\begin{prop}
The open orbit  in $\PP(\Delta_{14})$ is isomorphic with the homogeneous space 
$\spin_{14}/(G_2\times G_2)\rtimes \ZZ_2$.
\end{prop}

\proof It suffices to exhibit a transformation in $\spin_{14}$ that stabilizes $z_0$ and 
exchanges $V_7$ and $V'_7$. Remember from \cite{chevalley} 
that  $\spin_{14}$ embeds in the Clifford algebra of 
$V_{14}$ as the group generated by even products $g=v_1\cdots v_{2k}$ of norm one elements 
of $V_{14}$. Moreover the action on $V_{14}$ is obtained by mapping each $v_i$ to the 
corresponding orthogonal symmetry. Let $a_i=(e_i+f_i)/\sqrt{2}$ and $b_i=(e_i-f_i)/\sqrt{2}$
for $1\le i\le 7$, a set of vectors that constitute an orthonormal basis of $V_{14}$. Then
a straightforward computation shows that 
$$g=(a_1+a_4)(b_1+b_4)(a_2+a_5)(b_2+b_5)(a_3+a_6)(b_3+b_6)a_7b_7/8$$
belongs to the stabilizer of $z_0$ in $\spin_{14}$, and that its action on $V_{14}$ 
exchanges $V_7$ and $V'_7$. \qed 

\subsection{A multiplicative double point property}
Note that an orthogonal decomposition $V_{14}=V_7\oplus V'_7$ always determines a 
decomposition of $\Delta_{14}$ as $\Delta_7\otimes \Delta'_7$, as follows. By the 
Borel-Weil theorem, we can realize $\Delta_{14}$ as $H^0(S_{14},\cL_{14})$, where 
$\cL_{14}$ denotes the positive generator of the Picard group of the spinor variety
$S_{14}$. Similarly, we can realize $\Delta_{7}$ and $\Delta'_{7}$ as $H^0(S_{7},\cL_{7})$
and $H^0(S_{7},\cL'_{7})$, with similar notations. Points in $S_7$ and $S'_7$ are 
three-dimensional isotropic subspaces $E_3$ and $E'_3$ of $V_7$ and $V'_7$. Their
direct sum is still isotropic, and being of dimension six, it is contained in exactly 
two maximal isotropic spaces of $V_{14}$, one of each family. In particular this 
defines a regular map $\phi : S_7\times S'_7\rightarrow S_{14}$, such that $\phi^*\cL_{14}
=\cL_7\boxtimes\cL'_7$. By restriction this yields a map
$$\Delta_{14}\simeq H^0(S_{14},\cL_{14})\lra H^0(S_7\times S'_7, \phi^*\cL_{14})=
\Delta_7\otimes \Delta'_7.$$
This map is equivariant under $\spin(V_7)\times \spin(V'_7)$, and certainly non zero. 
Its target being irreducible, it has to be surjective, and then an isomorphism since
the source and target have the same dimension. 

We can summarize this discussion as follows. 

\begin{prop}\label{octooctonions}
Let $[\psi]$ be a generic element of $\PP(\Delta_{14})$. 
\begin{enumerate}
\item There exists a unique orthogonal decomposition $V_{14}=V_7\oplus V'_7$ preserved
by the stabilizer of $[\psi]$ in $\spin_{14}$.
\item Under the induced isomorphism $\Delta_{14}\simeq \Delta_7\otimes \Delta'_7$, 
we have $[\psi]=[\chi\otimes\chi']$ for $[\chi]$ and $[\chi']$ generic inside 
$\PP(\Delta_7)$ and $\PP(\Delta'_7)$.
%\item The stabilizer of $[\psi]$ is $(G\times G')\rtimes \ZZ_2$, where $G$ and $G'$ are
%the stabilizers of $[\chi]$ and $[\chi']$ in  $\spin(V_7)$ and $\spin(V'_7)$. In particular
%$G$ and $G'$ are both isomorphic to $G_2$. 
\end{enumerate}
\end{prop}

\noindent {\it Remarks.}
\begin{enumerate}
\item 
Note that $V_7$ must belong to the open subset of $G(7, V_{14})$ defined by the 
condition that the restriction of the quadratic form remains non degenerate. 
This open subset has dimension $49$, and for each choice of $V_7$ there are 
$7$ parameters for the generic $[\chi]$ and $[\chi']$. This yields the correct number 
$49+2\times 7=63$ of parameters for the open orbit in $\PP(\Delta_{14})$. 
Moreover, we get a remarkable partition of the open orbit in $\PP(\Delta_{14})$ by a $
49$-dimensional family of open subsets of $\PP^7\times\PP^7$.
\item The stabilizer of a generic point $[\chi]$  in $\PP\Delta_7$ is a copy of $G_2$, whose
action on $\Delta_7$ can be identified with the action of the latter on the Cayley algebra
$\OO$ (and $\chi$ becomes the unit in this algebra). As a consequence, a generic point in 
$\Delta_{14}$ allows to interprete it as the tensor product $\OO\otimes\OO$ 
of two Cayley algebras. It would be interesting to relate this observation to the work of 
Rosenfeld on the algebra of "octooctonions" \cite{rosenfeld}. 
\end{enumerate}
%\medskip
%By the last assertion of the proposition, the open orbit $\cO_0=\spin_{14}/(G\times G')\rtimes \ZZ_2$ in $\PP(\Delta_{14})$ has a natural double cover $\tilde\cO_0= Spin_{14}/G\times G'$.
%Choosing a point in the fiber of $\tau_0 : \tilde\cO_0\lra \cO_0$ over the point $[\psi]$ amounts to choosing $V_7$ or $V'_7$ as a prefered subspace in the direct sum $V_{14}=V_7\oplus V'_7$. In particular, $\tilde\cO_0$ is naturally endowed with two vector bundles $\cE_0$ and $\cE'_0$ 
%of rank seven, exchanged by the deck transformation $\iota$ of $\tau_0$. 

\subsection{A matrix factorization for the octic invariant}
We will construct later on a matrix factorization for the octic invariant $J_8$. 
A first explicit but cumbersome construction was obtained by Gyoja \cite{gyoja}.
Let us present a more direct approach. 

Our main observation is that, according to \cite{Lie}, the symmetric square of $\Delta_{14}$ contains a copy of $\wedge^3V_{14}$. This can be deduced from the Clifford action on the full 
spin representation, which decomposes into equivariant maps 
$$V_{14}\otimes \Delta_{14}\longrightarrow \Delta_{14}^\vee\quad \mathrm{and}\quad 
V_{14}\otimes \Delta_{14}^\vee \longrightarrow \Delta_{14}.$$
(Recall that $\Delta_{14}^\vee\simeq \wedge^-E$, on which $E\subset V_{14}$ acts by wedge product and  
$F\subset V_{14}$ by twice the contraction.) Composing those maps we get an equivariant morphism $\wedge^3V_{14}\otimes \Delta_{14}\hookrightarrow V_{14}\otimes V_{14}\otimes V_{14}\otimes \Delta_{14}\ra
V_{14}\otimes V_{14}\otimes \Delta_{14}^\vee\ra V_{14}\otimes \Delta_{14}\ra
\Delta_{14}^\vee$. Taking the transpose we get 
\begin{equation}
\Omega : \SS^2 \Delta_{14}\longrightarrow \wedge^3 V_{14}.
\end{equation}
To be more explicit, we can fix a basis $(v_1,\ldots ,v_{14})$ of $V_{14}$ (for example $(e_1, 
\ldots e_7, f_1,$ $ \ldots , f_7)$) and denote the dual basis by $(w_1,\ldots ,w_{14})$ 
(which would be $(f_1, \ldots f_7, e_1, $ $ \ldots , e_7)$ for the same example). Then 
\begin{equation}\label{omega}
\Omega_z=\sum_{i<j<k}\langle z,v_iv_jv_kz\rangle w_iw_jw_k \in \wedge^3 V_{14}.
\end{equation}
Here the natural pairing between $z\in\Delta_{14}=\wedge^+E$ and $z'\in\wedge^-E$ is defined
as the component of $z\wedge z'$ on $\wedge^7E$. 

\begin{lem}
The equivariant map $\Omega$ is non zero. 
\end{lem}

\proof We consider $\Omega$ as a quadratic form on $\Delta_{14}$, with values in 
$\wedge^3 V_{14}$, and evaluate it at a general point, that is, at the point $z_0$ 
of the open orbit. We get 
\begin{equation*}
\frac{1}{2}\Omega_{z_0} = e_{123}-e_{456}-f_{123} - f_{456} 
-\sum_{i=1}^6\epsilon_ie_{i7} f_i - \sum_{i=1}^6e_i f_{i7},
\end{equation*}
where $\epsilon_i=1$ for $i\le 3$ and $\epsilon_i=-1$ for $i\ge 4$. This is of course non zero. \qed

\medskip Observe moreover 
that  $\frac{1}{2}\Omega_{z_0}$  decomposes nicely as $\Omega-\Omega'$, 
where $\Omega\in\wedge^3V_7$ and $\Omega'\in\wedge^3V'_7$ are given by
\begin{equation*}
\Omega = e_{123} -f_{123}
+(\sum_{i=1}^3e_{i}\wedge f_i)\wedge (e_7-f_7), 
\end{equation*}
\begin{equation*}
\Omega' = e_{456} +f_{456} 
+(\sum_{i=4}^6e_{i}\wedge f_i)\wedge (e_7+f_7). 
\end{equation*}
Those forms $\Omega$ and $\Omega'$ are generic elements 
of $\wedge^3V_7$ and $\wedge^3V'_7$ (up to normalizations, they coincide with the generic
three-form explicited in \cite{cayley}). Recall that we recover $G_2$ as the stabilizer of 
such a generic form. 

\medskip
Our second ingredient will be the equivariant map
$$\Theta : \SS^2 (\wedge^3 V_{14})\longrightarrow End(V_{14})$$
obtained as follows. First embed $\wedge^3 V_{14}$ inside $V_{14}\otimes\wedge^2 V_{14}$,
Then recall that the quadratic form $q$ on $V_4$ induces a quadratic form 
$$q_2 : \SS^2 (\wedge^2 V_{14})\longrightarrow \CC,$$ 
whose polarization is given by the formula
$$q_2(u_1\wedge u_2,v_1\wedge v_2) 
=\det (q(u_i,v_j))_{1\le i,ij\le 2}.$$
Use this quadratic form  to define the composition
$$\Theta : \SS^2 (\wedge^3 V_{14})
\hookrightarrow  \SS^2 (V_{14}\otimes\wedge^2 V_{14})
\longrightarrow  \SS^2 V_{14}\otimes\SS^2(\wedge^2 V_{14})
\longrightarrow  \SS^2 V_{14}\hookrightarrow  End(V_{14}).$$
Finally, for $z\in\Delta$, let $M_z=\Theta (\Omega_z)\in End(V_{14})$. 
Using equation (\ref{omega}), we can compute explicitely
\begin{equation}
M_z = \sum_{k,\ell}\Big( \sum_{\substack{i<j, \\ i,j\ne k,\ell}}\langle z,v_iv_jv_kz\rangle \langle z,w_iw_jw_\ell z\rangle
\Big) w_kv_\ell.
\end{equation}
%\medskip 

\begin{theo}\label{mf14}
The pair $(M,M)$
is a matrix factorization of the octic invariant  $J_8$ of $\Delta_{14}$. 
\end{theo}

\proof We just need to check that $M_{z_0}^2$ is a non zero multiple of the identity. So 
let us compute $M_{z_0}$. We have seen that  $\frac{1}{2}\Omega_{z_0}=\Omega-\Omega'$, where 
$\Omega$ and $\Omega'$ belong to $\wedge^3V_7$ and  $\wedge^3V'_7$, respectively. 

\begin{lem}
$\Theta (\Omega-\Omega')=\Theta (\Omega)+\Theta (\Omega')$.
\end{lem}

\proof Indeed, let $v_1, \ldots , v_7$ and $v'_1, \ldots , v'_7$ be basis of $V_7$ and $V'_7$,
respectively. The polarisations of $\Omega$ and $\Omega'$ in $V_{14}\otimes \wedge^2V_{14}$
will respectively be of the form $\sum_{i=1}^7v_i\otimes \omega_i$ and $\sum_{i=1}^7v'_i\otimes \omega'_i$, for some tensors $\omega_i\in\wedge^2V_7$ and $\omega'_i\in\wedge^2V'_7$. When
we apply $\Theta$ and take the image in $S^2V_{14}$, the mixed terms are of the form $q_2(\omega_i,\omega'_j)v_iv'_j$. But $q_2(\omega_i,\omega'_j)$ is always
zero since $V_7$ and $V'_7$ are orthogonal. \qed 

\medskip In order to compute $\Theta (\Omega)$, we first send
$\Omega$ to $V_{14}\otimes\wedge^2 V_{14}$ by polarizing it. Let $e_0=e_7-f_7$. We get 
$$\begin{array}{rcl} 
\bar\Omega &=& e_1\otimes (e_{23}+f_1e_0)+e_2\otimes (e_{31}+f_2e_0)+e_3\otimes (e_{12}+f_3e_0)+\\
 & & \hspace{1cm} +f_1\otimes (f_{23}-e_1e_0)+ f_2\otimes (f_{31}-e_2e_0)+f_3\otimes (f_{12}-e_3e_0)+ \\
 & & \hspace{2cm} +e_0\otimes (e_1f_1+e_2f_2+e_3f_3).
\end{array}$$
Now recall that $q(e_i,f_i)=1$ for all $i$, while $q(e_0)=-2$; moreover $q$ evaluates
to zero on any other pair of basis vectors. We deduce that $q_2(e_1f_1+e_2f_2+e_3f_3)=3$,
$$q_2(e_{23}+f_1e_0,f_{23}-e_1e_0)=q_2(e_{31}+f_2e_0,f_{31}-e_2e_0)=q_2(e_{12}+f_3e_0,f_{12}-e_3e_0)=3,$$
and that all the other scalar products are zero. This yields
$$\Theta (\Omega)=3e_0^2-6e_1f_1-6e_2f_2-6e_3f_3 \in S^2V_{14}.$$
With respect to the quadratic form $q$, the dual basis of $(e_0,e_1,e_2,e_3,f_1,f_2,f_3)$ is 
$(-\frac{1}{2}e_0,f_1,f_2,f_3,e_1,e_2,e_3)$. Considered as an element of $End(V_{14})$, 
the tensor $\Theta (\Omega)$ is thus exactly $-6\pi_{V_7}$, where $\pi_{V_7}$ denotes 
the orthogonal projection to 
$V_7$. A similar computation shows that $\Theta (\Omega')$
is  $+6\pi_{V'_7}$. We  finally get
$$M_{z_0}=24(\pi_{V'_7}-\pi_{V_7}),$$
whose square is $576$ times  the identity. This concludes the proof. \qed

\medskip\noindent {\it Remarks}.
\begin{enumerate}
\item Once we have observed that $\Theta (\Omega_{z_0})=\Theta (\Omega)+\Theta (\Omega')$,
we can in fact conclude without any extra computation. Indeed, $\Theta (\Omega)$ is an element
of $S^2V_7$ that must be preserved by the stabilizer of $\Omega$, hence by a copy of $G_2$. 
But up to scalar there is a unique such element. Moreover we already know one: the restriction
to $V_7$ of the quadratic form on $V_{14}$. The same being true for $\Omega'$, we conclude 
that there exist scalars $a$ and $a'$ such that 
$$M_{z_0}=a\pi_{V_7}+a'\pi_{V'_7}.$$
But the trace of $M_{z_0}$ must be zero, otherwise we would get a non trivial quartic
invariant on $\Delta_{14}$, and we know there is none. So $a+a'=0$, and the square 
of $M_{z_0}$ is a homothety. 
\item
Let us also compute $M_{z_1}$. We start by computing $\Omega_{z_1}$: 
$$\begin{array}{rcl}
\frac{1}{2}\Omega_{z_1} &= & e_{123}+e_{156}+e_{246}-f_{135}-f_{234}-f_{456}+ \\
 & & \hspace*{2cm} +(e_2f_5-e_1f_4+e_3f_6)e_7-(\sum_{i=1}^6e_if_i)f_7.
 \end{array}$$
 Polarizing, we get the following tensor $\Theta_1$ in $V_{14}\otimes\wedge^2V_{14}$:
 $$\begin{array}{r}
 e_1\otimes (e_{23}+e_{56}+f_4e_7-f_{17})+e_2\otimes (e_{46}-e_{13}-f_5e_7-f_{27})\\
+e_3\otimes (e_{12}-f_6e_7-f_{37}) 
-e_4\otimes (e_{26}+f_{47})-e_5\otimes (e_{16}+f_{57}) \\+e_6\otimes (e_{15}+e_{24}-f_{67})+ 
e_7\otimes (e_2f_5-e_1f_4  +e_3f_6)  \\+f_1\otimes (e_1f_7-f_{35})+f_2\otimes (e_2f_7-f_{34})+
f_3\otimes (e_3f_7+f_{15}+f_{24}) \\ + f_4\otimes (e_{17}+e_4f_7-f_{23}-f_{56})+ 
-f_5\otimes (e_{27}-e_5f_7+f_{13}-f_{46}) \\ -f_6\otimes (e_{37}-e_6f_7+f_{45})-f_7\otimes (\sum_{i=1}^6e_if_i).
 \end{array}$$
If we write $\Theta_1=\sum_{i=1}^7(e_i\otimes a_i+f_i\otimes b_i)$, it is straightforward to check that 
the only non zero scalar products between the two-forms $a_i,b_j$ are the following:
$$q_2(a_1,b_4)=-2, \; q_2(a_2,b_5)=2, \; q_2(a_6,b_3)=-2, \; q_2(b_7)=6. $$
We thus finally get $M_{z_1}$ as the following element of $\SS^2V_{14}$:
$$M_{z_1}=8(3f_7^2-e_1f_4+e_2f_5+e_6f_3).$$
As an endomorphism of $V_{14}$, $M_{z_1}$ has for image and kernel the same vector space 
$V_7=\langle f_7,e_1,f_4,e_2,f_5,e_6,f_3\rangle$. In particular, the square of $M_{z_1}$ is
zero, in agreement with the fact that $J_8(z_1)=0$. Note that $V_7$
is isotropic; moreover, since it meets $E$ in odd dimension, it belongs to the same family of maximal
isotropic subspaces, which is embedded in the other projectivized half-spin representation $\PP (\Delta_{14}^\vee)$. 
This is in agreement with the fact that the octic invariant hypersurface in $\PP (\Delta_{14})$ can be 
obtained as the projective dual variety of the closed orbit $S'_{14}\subset \PP (\Delta_{14}^\vee)$. In particular,
the open orbit inside the octic is naturally fibered over $S'_{14}$, and our $z_1$ must be sent to 
$V_7$ by this fibration. 
\end{enumerate}

\smallskip
Let us summarize what we have proved so far, which 
is amazingly similar to what happens for $S_{12}$ and $\wedge^3V_6$, see \cite[sections 3.3 and 3.4]{ulrich}. 

\begin{prop}
\begin{enumerate} 
\item Let $[z]$ belong to the open orbit in $\PP (\Delta_{14})$, and let $(V_7,V'_7)$ be the associated
pair of orthogonal non degenerate subspaces of $V_{14}$. Then the associated three-form 
$\Omega_z$ is the sum of generic three-forms $\Omega\in\wedge^3V_7$ and $\Omega'\in\wedge^3V'_7$. 
Moreover 
$$M_z=m_z(\pi_{V'_7}-\pi_{V_7}), \qquad with \quad J_8(z)=m_z^2.$$
\item  Let $[z]$ belong to the open orbit in the invariant hypersurface $(J_8=0)$ of $ \PP (\Delta_{14})$, 
and let $V_7$ be the associated maximal isotropic subspace of $V_{14}$. Then $\Omega_z$ belongs to $\wedge^2V_7\wedge V_{14}$, 
and $M_z$ is a square zero endomorphism whose image and kernel are both equal to $V_7$. 
\end{enumerate}
\end{prop}

Beware that in $(1)$, the pair $(V_7, V'_7)$ is not ordered. Permuting $V_7$ and $V'_7$ 
changes the sign of $m_z$, so only its square is well-defined. And $J_8$ is not
globally a square. 

\subsection{Double octics}
Consider the double cover  $\pi : D\lra \PP(\Delta_{14})$, branched over the octic hypersurface $(J_8=0)$. This double cover can be interpreted as the octic hypersurface $J_8(z)-y^2=0$ in the weighted projective space 
$\tilde{\PP}=\PP(1^{64},4)$. Moreover it inherits an action of the spin group $\spin_{14}$. 

\begin{prop} 
\begin{enumerate}
\item The double cover $D$ is smooth in codimension $5$. 
\item Its smooth locus $D_0$ is endowed with two rank seven equivariant vector bundles $\cE_+$ and $\cE_-$, exchanged by the 
deck transformation $\iota$ of the double cover, which are two orthogonal subbundles of the trivial bundle $\cV$ 
with fiber $V_{14}$. Moreover there are exact sequences of vector bundles 
$$0\lra \cE_\pm\lra \cV\lra \cE_\mp ^\vee \lra 0.$$
\end{enumerate}
\end{prop}

The matrix factorization $(M_z,M_z)$ of $J_8(z)$ upgrades to a matrix factorization $(M_z+yId,M_z-yId)$ of 
$J_8(z)-y^2$. The bundles $\cE_+$ and $\cE_-$ can be defined at the point $[z,y]$ as the kernels of 
$M_z+yId$ and $M_z-yId$, respectively. Moreover, denote by $\tilde{\PP}_0$ the complement in $\tilde{\PP}$ of
the singular locus of the hypersurface $(J_8=0)$. Denote by $j_0$ the embedding of $D_0$ inside 
$\tilde{\PP}_0$. Then we have the following exact sequences of sheaves on $\tilde{\PP}_0$:
$$0\lra V_{14}\otimes \cO_{\tilde{\PP}_0}(-4) \stackrel{M_z\pm yId}{\lra} V_{14}\otimes \cO_{\tilde{\PP}_0}
\lra j_{0*}\cE_\mp ^\vee\lra 0. $$

\medskip
If $L\subset\PP(\Delta_{14})$ is a general linear subspace of dimension at most four, then it is contained
in $D_0$. Over $L$ we then get a double cover $X_L$ with two rank seven vector bundles $\cE_L$ and $\cE'_L$. 

\medskip\noindent {\bf Conjecture}. {\it For a general three dimensional subspace $L\subset\PP(\Delta_{14})$, 
the vector bundles $\cE_L$ and $\cE'_L$ over the double cover $X_L$, are spherical.}

\medskip
Note that in this case, $X_L$ is a Calabi-Yau threefold.
By spherical, one means that the bundles of endomorphisms of $\cE_L$ and $\cE'_L$ has
the same cohomology as a three-dimensional sphere:
$$H^q(X,End(\cE_L))=H^q(X,End(\cE'_L))=\delta_{q,0}\CC\oplus \delta_{q,3}\CC.$$
As already discussed in \cite{im}, one can expect that a general
double octic threefold $X$ can always be represented as a section $X_L$ of the octic in $\PP(\Delta_{14})$. 
Moreover there should exist only a finite number $N$ of such representations (up to isomorphism).
This is the exact analogue of the statement, due to Beauville and Schreyer, that a general quintic threefold can be represented as a Pfaffian, in only finitely many ways \cite[Proposition 8.9]{beauville}.

Moreover, our construction 
of the supposedly spherical vector bundles  $\cE_L$ and $\cE'_L$ would parallel those of 
spherical bundles of rank seven on the generic cubic sevenfold \cite{im}, and of rank six 
on the generic double quartic fivefold \cite{abuaf}.

Recall that a spherical object in the derived category of coherent sheaves of an algebraic variety
defines a non trivial auto-equivalence of this category, called a spherical twist \cite{st}. If the previous conjecture is true, 
we would thus get, on a general octic threefold, $N$ pairs of spherical vector bundles, generating a certain group of 
symmetries of the derived category. It would be interesting to determine $N$, and this symmetry group. 

The reconstruction problem also looks very intringuing: starting from a general octic threefold $X$, its branch divisor $D$, and the vector 
bundle $\cE_L$, how can we reconstruct the embedding of $D$ as a linear section of the invariant octic in $\PP(\Delta_{14})$?

\section{Relations with $\ZZ$-gradings of Lie algebras}

\subsection{Morphisms from gradings}
Clerc was the first to realize that $\ZZ$-gradings can be useful to determine 
certain invariants \cite{clerc}. He starts with a simple Lie algebra $\fg$ whose
adjoint representation is fundamental. In other words (once we have fixed a Cartan
and a Borel subalgebra), the highest root $\psi$ is a fundamental weight. The 
corresponding simple root $\alpha$ defines a $\ZZ$-grading on $\fg$, 
of length five:
$$\fg = \fg_{-2}\oplus  \fg_{-1}\oplus  \fg_{0}\oplus  \fg_{1}\oplus  \fg_{2}.$$
In this grading, $\fg_2=\fg_\psi$ and $\fg_{-2}=\fg_{-\psi}$ are one dimensional. 
Once we have chosen generators $X_\psi$ and $X_{-\psi}$, we can define a $\fg_0$-equivariant 
polynomial function $J_4$ on $\fg_1$, homogeneous of degree four, by the relation 
$$ad(z)^4X_{-\psi} = J_4(z)X_\psi, \qquad z\in\fg_1.$$
For $\fg$ exceptional, $J_4$ generates the space of semi-invariants for the action of $\fg_0$ on 
$\fg_1$. (Recall that $\fg_0$ is not semisimple but only reductive, with a non trivial center. 
That $J_4$ is semi-invariant means that $\fg_0$ acts on it only through multiplication by some
character.)

\smallskip
There exist other gradings of length five of simple Lie algebras, notably of the 
exceptional ones, such that $\fg_2$ has  dimension bigger than one. In this case,
the very same idea yields morphisms 
$$Sym^4\fg_1\lra Hom(\fg_{-2},\fg_2) \quad \mathrm{and}\quad 
Sym^4\fg_{-1}\lra Hom(\fg_{2},\fg_{-2}).$$
The representations $\fg_1$ and $\fg_{-1}$, as well as $\fg_2$ and $\fg_{-2}$, 
are in perfect duality through the Cartan-Killing form on $\fg$. But it is often
the case that $\fg_2$ is in fact self-dual. The quite unexpected fact, discovered 
in \cite{abuaf}, is that we can construct matrix factorizations from the resulting maps.

\subsection{A magic square of matrix factorizations}
A remarkable framework in which we will obtain matrix factorizations is that of the 
Tits-Freudenthal magic square (see \cite[Section 4.3]{baez} and references therein). 
Either in its algebraic, or in its geometric versions, 
this magic square (and its enriched triangular version due to Deligne and Gross \cite{dg}) 
encodes all sorts of surprising phenomena related to the exceptional Lie algebras. 
We will describe another one in this section. 

Remember  that the magic square has its lines and
columns indexed by the normed algebras $\RR, \CC, \HH, \OO$ and is symmetric at the 
algebraic level, in the sense that there is a way to associate to a pair $(\AA,\BB)$ 
a semisimple Lie algebra $\fg(\AA,\BB)$, with $\fg(\AA,\BB)=\fg(\AA,\BB)$. Here is the magic square
over the complex numbers:

\begin{center}\begin{tabular}{cccc}
$\fso_3$ & $\fsl_3$ & $\fsp_6$ &
$\ff_4$\\ $\fsl_3$ & $\fsl_3\times \fsl_3$ &
$\fsl_6$ & $\fe_6$ \\ $\fsp_6$ & $\fsl_6$ & $\fso_{12}$ & $\fe_7$ 
\\ $\ff_4$ & $\fe_6$ & $\fe_7$ & $\fe_8$ \end{tabular}\end{center}

Freudenthal enhanced this construction by associating to each pair $(\AA,\BB)$
some special geometry, in a way which is not symmetric, but uniform along the lines
(see \cite[Section 2.1]{LMfreud} and references therein).
This means that for each $\AA$, one has four special geometries associated to 
the pairs $(\AA,\RR)$, $(\AA,\CC)$, $(\AA,\HH)$, $(\AA,\OO)$, with completely
similar formal properties. Each of these geometries is made of elements that we call
F-points (for the first line of the magic square), plus F-lines (for the second 
line), plus F-planes (for the third line), plus F-symplecta (for the fourth line).
Moreover, each type of elements involved in the special geometry of the pair $(\AA,\BB)$  
is parametrized by a rational homogeneous space $G/P$, where $G$ is an algebraic group 
whose Lie algebra is precisely $\fg(\AA,\BB)$. The parabolic subgroup $P$,
usually maximal, depends on the type of the element. All these data can be encoded in
the following diagrams:

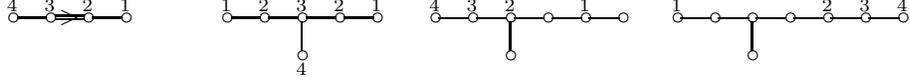
\begin{center}\begin{tabular}{cccc}
\setlength{\unitlength}{2.5mm}
\begin{picture}(10,3)(-.5,-1)
\multiput(0,0)(2,0){4}{$\circ$}
\put(0.55,.35){\line(1,0){1.55}}
\put(4.55,.35){\line(1,0){1.55}}
\multiput(2.55,.25)(0,.2){2}{\line(1,0){1.55}} \put(2.75,0){$>$}
\put(0,.7){$\scriptstyle{4}$}
\put(2,.7){$\scriptstyle{3}$}
\put(4,.7){$\scriptstyle{2}$}
\put(6,.7){$\scriptstyle{1}$}
\end{picture} &
\setlength{\unitlength}{2.5mm}
\begin{picture}(10,3)(0,-1)
\multiput(0,0)(2,0){5}{$\circ$}
\multiput(0.55,.35)(2,0){4}{\line(1,0){1.55}}
\put(4.3,-1.4){\line(0,1){1.5}}
\put(4,-2){$\circ$}
\put(0,.7){$\scriptstyle{1}$}
\put(2,.7){$\scriptstyle{2}$}
\put(4,.7){$\scriptstyle{3}$}
\put(6,.7){$\scriptstyle{2}$}
\put(8,.7){$\scriptstyle{1}$}
\put(4,-2.7){$\scriptstyle{4}$}
\end{picture} &
\setlength{\unitlength}{2.5mm}
\begin{picture}(11,3)(0.8,-1)
\multiput(0,0)(2,0){6}{$\circ$}
\multiput(0.5,.3)(2,0){5}{\line(1,0){1.6}} \put(4,-2){$\circ$}
\put(4.3,-1.4){\line(0,1){1.5}}
\put(0,.7){$\scriptstyle{4}$}
\put(2,.7){$\scriptstyle{3}$}
\put(4,.7){$\scriptstyle{2}$}
\put(8,.7){$\scriptstyle{1}$}
\end{picture} &
\setlength{\unitlength}{2.5mm}
\begin{picture}(11,3)(0.8,-1)
\multiput(0,0)(2,0){7}{$\circ$}
\multiput(0.5,.3)(2,0){6}{\line(1,0){1.6}} \put(4,-2){$\circ$}
\put(4.3,-1.4){\line(0,1){1.5}}
\put(0,.7){$\scriptstyle{1}$}
\put(10,.7){$\scriptstyle{3}$}
\put(8,.7){$\scriptstyle{2}$}
\put(12,.7){$\scriptstyle{4}$}
\end{picture} \\
\end{tabular}\end{center}\medskip
This must be interpreted as follows. Each of these diagrams corresponds to one column in the 
magic square. Recall that a Dynkin diagram encodes a semisimple Lie algebra
(or Lie group, up to finite covers), and that a subset of vertices encodes a conjugacy class of 
parabolic subgroups.  The geometry associated to the box $(i,j)$, on the $i$-th line and the $j$-th column, 
is defined by considering the $j$-th Dynkin diagram above and suppressing the vertices numbered by
integers bigger than $i$. This gives a Dynkin diagram $D_{i,j}$ marked by integers from $1$ to $i$. 
The elements of the corresponding geometry are then parametrized by the homogeneous spaces $G/P(k)$,
$1\le k\le i$, where $G$ is a semisimple Lie group with Dynkin diagram $D_{i,j}$, and $P(k)$ is a 
parabolic subgroup defined, up to conjugacy, by the vertices of $D_{i,j}$ marked by $k$. The spaces
$G/P(1)$ parametrize F-points, while the $G/P(2)$ parametrize F-lines (for $i\ge 2$), the 
$G/P(3)$ parametrize F-planes (for $i\ge 3$), and the $G/P(4)$ parametrize F-symplecta
(for $i=4$). For example, the Dynkin diagram $D(3,2)$ has type $A_5$, with the marks $1$ at its extremities; 
so that the space of F-points for $\AA=\HH$ and $\BB=\CC$ is $A_5/P_{1,5}$, the flag variety of incident
points and hyperplanes in $\PP^5$.

Let us focus on the square of homogeneous spaces $X(\AA,\BB)$ parametrizing the 
F-points of the Freudenthal geometries, namely:

\begin{center}\begin{tabular}{cccc} $A_1/P_1$ & $A_2/P_{1,2}$ &
$C_3/P_2$ & $F_4/P_1$  \\ $A_2/P_1$ &
$A_1\times A_1/P_{1,1}$ & $A_5/P_2$ & 
$E_6/P_1$ \\
$C_3/P_1$ & $A_5/P_{1,5}$ & $D_6/P_2$ & $E_7/P_1$  \\ 
$F_4/P_1$ &
$E_6/P_{1,6}$ & $E_7/P_6$ & $E_8/P_1$  \end{tabular} \end{center} 
Now, we define a five-step grading of $\fg(\AA,\BB)$ as follows. As we have just explained,
each $X(\AA,\BB)$ is $G/P_I$, where the standard parabolic subgroup $P_I$ of $G$ 
is defined by the subset $I$ of the set $\Delta$ of simple roots (which correspond bijectively 
with the vertices of the Dynkin diagram). Let 
$\omega_I^\vee=\sum_{i\in I}\omega_i^\vee$ denote the sum of the corresponding 
fundamental coweights. If we express a root $\alpha$ as a linear combination of simple roots, 
$\alpha = \sum_{j\in\Delta}n_j\alpha_j$, then $\omega_I^\vee(\alpha)=\sum_{i\in I}n_i$. 
The associated $\ZZ$-grading of $\fg=\fg(\AA,\BB)$  is 
$$\fg_k= \delta_{k,0}\ft \oplus (\bigoplus_{\omega_I^\vee(\alpha)=k}\fg_\alpha ). $$
%In the two squares below we provide the  type of $\fg_0$ (which is always reductive) 
%and the representations $\fg_1,\fg_2$. 
It turns out that $\fg_0$ is always made of orthogonal
Lie algebras, that $\fg_1$ is always a spin representation, while $\fg_2$ is a natural
representation, in particular self dual. Moreover $\fg_k=0$ for $k>2$. 
It was argued in \cite[Proposition 3.2]{LMfreud}
that $\fg_1$ should be thought of as $\AA\otimes\BB$ and $\fg_2$ as $\AA_0\oplus\BB_0$,
where $\AA_0$ is the hyperplane of imaginary elements in $\AA$. Here is the table giving 
the semisimple part of $\fg_0$:

\begin{center}\begin{tabular}{cccc}
$0$ & $0$ & $\fsl_2^2$ & $\fspin_7$\\ 
$\fsl_2$ & $\fsl_2^2$ &
$\fsl_2\times\fsl_4$ & $\fspin_{10}$ \\ 
$\fsp_4$ & $\fsl_4$ & $\fsl_2^2\times\fsl_3$ & $\fspin_{12}$ \\ 
$\fspin_7$ & $\fspin_8$ & $\fsl_2\times\fspin_{10}$ & $\fspin_{14}$ 
\end{tabular}\end{center}

%On the fourth line, we choose the spaces that parametrize F-points. On the third 
%line, we take those that parametrize F-lines. On the second line, we consider the 
%spaces that parametrize F-flags, that is, pairs of incident points and lines. It
%turns out that the resulting rectangle is partially symmetric, in the sense that 
%$X(\AA,\BB)=X(\BB,\AA)$ for $\AA,\BB\ne\RR$. So we complete it in a symmetric way,
%and we let $X(\RR,\RR)=?$. We get the following square:

\medskip

\begin{theo}\label{magicmf}
For each pair $(\AA,\BB)$ the representations  $\fg_1,\fg_2$ of the reductive Lie algebra 
$\fg_0$ have the following properties:
\begin{enumerate}
\item $\fg_1$ is prehomogeneous under the action of $\fg_0$, and the generic stabilizer is 
$aut(\AA)\times aut(\BB)$,
%(possibly up to $\ZZ_2$ when $\AA=\BB$),
\item there exist equivariant morphisms 
$$P,Q : Sym^4\fg_1\lra Sym^2\fg_2\hookrightarrow End(\fg_2)$$ 
such that $(P,Q)$ is a matrix factorization of a semi-invariant 
$J_8$ on $\fg_1$. 
\end{enumerate}
\end{theo}
 
This statement generalizes Theorem \ref{mf14}, which corresponds to the pair $(\OO,\OO)$. 
All the other cases are somewhat degenerate, in the following ways. 
\begin{enumerate}
\item For the pairs $(\OO,\CC)$ and $(\OO,\HH)$ of the last line, the fundamental semi-invariant has degree four, and we get a matrix factorization of its square. 
\item For the pair $(\HH,\OO)$ of the third line, $\fg_2$ is one dimensional, the fundamental
invariant has degree four and admits a matrix factorization induced by the 5-grading, as 
in \cite{abuaf}.
\item For the pair $(\CC,\OO)$ of the second line, there is no non trivial semi-invariant 
and $\fg_2$ is actually zero.
\end{enumerate}
The matrices $P_z$ and $Q_z$ will always be obtained as $ad(z)^4$, up to some homothety, 
but the proof of the Theorem is unfortunately a case by case check. In the next section we 
will discuss the first degenerate case, that of $V_2\otimes\Delta_{10}$ where $V_2$ is two-dimensional. 
The other cases are similar and simpler. 

\subsection{The case of $\CC^2\otimes\Delta_{10}$}
Here the action of $GL(V_2)\times \spin_{10}$ is
prehomogeneous, and there is a quartic semi-invariant $J_4$ \cite{kwE7}. Note that this is
the representation used by Hitchin in order to defined (a substitute for) $GL_2(\OO)$, with
$J_4$ playing the r\^ole of the determinant \cite{hitchin}. 

An important observation 
is that the natural map, again induced by Clifford multiplication, yields an isomorphism
$$\wedge^2\Delta_{10}\simeq \wedge^3V_{10}.$$
We therefore get a morphism of $SL(V_2)\times \spin_{10}$-modules
\begin{equation}
\Omega : S^2(V_2\otimes\Delta_{10})\lra \wedge^2V_2\otimes\wedge^2\Delta_{10}\simeq \wedge^3V_{10}.
\end{equation}
This in turn induces a morphism 
\begin{equation}
M : S^4(V_2\otimes\Delta_{10})\lra S^2(\wedge^3V_{10})\lra S^2V_{10}\hookrightarrow 
End(V_{10}),
\end{equation}
where the first arrow is given by the square of $\Omega$, and the second arrow is induced 
by the quadratic form on $V_{10}$. 
Let us compute this morphism explicitely. According to \cite{kwE7}, a generic element of 
$V_2\otimes\Delta_{10}$ is given by 
$$z= v_1\otimes (1+e_{1234})+v_2\otimes (e_{1235}+e_{45}).$$
Letting $e_0=e_4+f_4$ and $f_0=f_4-e_4$, the associated three-form is 
$$\Omega_z = e_{123}+f_{123}+e_0(e_1f_1+e_2f_2+e_3f_3)+f_0e_5f_5.$$
Observe that this decomposes as the sum of $\Omega\in  \wedge^3V_7$ and $\Omega'=f_0e_5f_5
\in  \wedge^3V_3$, where $V_7=\langle e_0,e_1,e_2,e_3,f_1,f_2,f_3\rangle$ and $V_3=\langle
f_0,e_5,f_5\rangle$ are orthogonal spaces, on which the quadratic form $q$ is non degenerate. 
The associated element in $S^2V_{10}$ is $6e_1f_1+6e_2f_2+6e_3f_3+3e_0^2+f_0^2-2e_5f_5$, 
and since $q(e_0)=2$ and $q(f_0)=-2$, the corresponding endomorphism is 
$$M_z=6\pi_{V_7}-2\pi_{V_3}\in End(V_{10}),$$
where $\pi_{V_7}$ and $\pi_{V_3}$ are the two projections relative to the direct sum 
decomposition $V_{10}=V_7\oplus V_3$. This endomorphism has non zero-trace, which allows to define a
degree four invariant $J_4(z):=\mathrm{trace} (M_z)$. We then get 
$$(M_z-\frac{1}{18}J_4(z)Id_{V_{10}})^2=\frac{1}{162}J_4(z)^2Id_{V_{10}},$$
a matrix factorization of the octic $J_8=\frac{1}{162}J_4^2$. 

%\medskip
%Note that $ad(z)^4\in Hom(\fg_{-2},\fg_2)$ is invariant under the stabilizer of $z$, 
%which we know to be
%isomorphic to $\fsl_2\times\fg_2$. Moreover, this stabilizer acts on $V_{10}=g_{-2}=g_2$ through
%the direct sum 
%decomposition $V_{10}=V_7\oplus V_3$, the $\fg_2$-factor acting trivially on $V_3$ and irreducibly
%on $V_7$, and conversely for the $\fsl_2$-factor (see \cite[(5.42)]{sk}. By the Schur lemma, there 
%must therefore exist scalars $a_7$ and $a_3$, independent of $z$, such that 
%$$ad(z)^4=J_4(z)(a_7\pi_{V_7}+a_3\pi_{V_3})\in End(V_{10})\simeq Hom(g_{-2},g_2).$$
%As a consequence, there must exist a scalar $a_0$ such that $M_z= ad(z)^4+a_0J_4(z)Id_{V_{10}}$. 
%Determining $a_0$ would require to dig deeper into the structure of $\fe_7$, and we will
%not do that. 

\subsection{Three-forms in seven variables}
A sporadic case to which the same circle of ideas can be applied is that of the 
degree seven invariant of $\wedge^3V_7$. This heptic invariant has been known for a long 
time, as given by the equation of the projective dual of the Grassmannian $G(3,V_7)$, 
or of the complement of the open $GL_7$-orbit of $\wedge^3V_7$
consisting of forms with stabilizer bigger than $G_2$. This case is associated
with the grading of length five of $\fe_7$ defined by the simple root $\alpha_2$:
$$\fe_7 = V_7\oplus \wedge^3V_7^\vee \oplus \fgl(V_7)\oplus \wedge^3V_7\oplus V_7^\vee.$$
(For simplicity we wrote this grading as a decomposition into $\fsl(V_7)$-modules.)
For any $\omega \in \wedge^3 \V$ and $y \in \V^\vee$, we denote by $i_y(\omega) \in \wedge^2 \V$ the contraction of $\omega$ with $y$. Then $P(\omega,y)=\frac{1}{6}\omega \wedge i_y(\omega) \wedge i_y(\omega)$ belongs to $\wedge^7V_7\simeq \CC$. This defines an equivariant morphism
\begin{equation}
P :   \SS^3 (\wedge^3 \V) \longrightarrow \SS^2 \V \otimes \mathrm{det} \V .
\end{equation}

On the other hand, recall that $\wedge^3 \V\simeq (\wedge^4 \V)^\vee$ and that $(\wedge^4 \V)^\vee$ is a submodule of $\V ^\vee\otimes (\wedge^3\V)^\vee$. Moreover the natural map $End(\V)\longrightarrow End(\wedge^3 \V)$ has for transpose an equivariant map $End(\wedge^3 \V)=(\wedge^3\V)^\vee\otimes (\wedge^3\V)\longrightarrow End(\V)$.
This yields a morphism
\begin{equation}
\theta :   \SS^2 (\wedge^3 \V) \longrightarrow \V^\vee \otimes End(\V)\otimes \mathrm{det} \V .
\end{equation}
Taking its symmetric square and composing with the trace map, we get 
\begin{equation}
R :   \SS^4 (\wedge^3 \V) \longrightarrow \SS^2\V^\vee \otimes \SS^2 End(\V)\otimes (\mathrm{det} \V )^2 
\longrightarrow \SS^2\V^\vee \otimes (\mathrm{det} \V )^2 .
\end{equation}

\medskip\noindent {\it Remark.} 
The morphism $R$ is induced from the quartic $\mathrm{SL}_6$-invariant of 
$\wedge^3 U_6$, in the following way. This quartic corresponds to a 
$\mathrm{GL}_6$-equivariant linear map:
\begin{equation*}
\psi_{\U} : \SS^4 (\wedge^3 \U) \longrightarrow \mathrm{det}(\U)^{\otimes 2}.
\end{equation*}
Now the Borel-Weil Theorem implies that $\wedge^3 \V = H^0(\mathbb{P}(\V), \wedge^3 Q)$, 
where $Q$ is the tautological quotient on $\mathbb{P}(\V)$, the projective space of lines 
in $\V$. We can then consider the composition 
\begin{equation*}
\SS^4  H^0(\mathbb{P}(\V), \wedge^3 Q) \longrightarrow H^0(\mathbb{P}(\V), \SS^4 (\wedge^3 Q)) \longrightarrow H^0(\mathbb{P}(\V), \mathrm{det}(Q)^{\otimes 2}), 
\end{equation*}
where the first arrow is surjective because $Q$ is globally generated, and the second
arrow is defined by $\psi_Q$. 
Finally, the Borel-Weil Theorem implies that 
\begin{equation*}
H^0(\mathbb{P}(\V), \mathrm{det}(Q)^{\otimes 2}) = \SS^2 \V^\vee \otimes \det (\V)^{\otimes 2},
\end{equation*}
and the resulting morphism is nothing else but $R$. 

\medskip
As was first done by Gyoja and Kimura \cite{kimura}, 
we can now define the heptic semi-invariant $J_7$ on
$\wedge^3V_7$ by contracting $R$ with $P$:
$$J_7(z)=\langle P_z, R_z\rangle \in (\det\V)^{\otimes 3}.$$
But in fact a much stronger statement is true. Observe that $\SS^2 \V $ is a submodule
of $Hom(\V^\vee, \V)$, as well as $\SS^2 \V^\vee $ is a submodule
of $Hom(\V, \V^\vee)$. In other words, we can consider $R_z$ and  $P_z$
as symmetric morphisms from $\V$ to $\V^\vee$ and from $\V^\vee$ to $\V$, respectively. 
The following result appears in \cite[Example 2.6]{kimurabook}. We give a short proof,
without computation. 

\begin{theo} \label{hepticmf}
The pair of symmetric morphisms $(P, R)$ is a matrix factorization
of the heptic semi-invariant $J_7$ of $\wedge^3 \V$.
\end{theo}

\proof Because of the quasi-homogeneity, it is enough to check this assertion at a generic
point $z$ of $\wedge^3\V$. The stabilizer of this form is then a copy of $G_2$. The quadratic
forms $R_z$ and  $P_z$ on $\V ^\vee$ and $\V$ must be preserved by 
this stabilizer. But up to scalar, there is a unique quadratic form on $\V$ (or its dual) 
preserved by $G_2$. After identifying $\V$ with its dual through this non degenerate quadratic
form (which depends on $z$), we get $R_z$ and  $P_z$ as (non zero) homotheties, and our statement immediately follows. \qed

\medskip\noindent {\it Remark.} 
In many respects, the case of $(\spin_{14},\Delta_{14})$ is a doubled version of the case of $(\sl_7,\wedge^3\CC^7)$. 
Similarly, $(\spin_{12},\Delta_{12})$ is a doubled version of $(\sl_6,\wedge^3\CC^6)$ (they both appear on the 
same line of Freudenthal's magic square, see \cite{LMfreud}), and 
  $(\spin_{10},\Delta_{10})$ is a doubled version of $(\sl_5,\wedge^2\CC^5)$ (see
\cite[Introduction]{doublespinor} for more details). 
The second pair has a clear complex-quaternionic interpretation. It would be very nice to find a similar, or may be 
quaternio-octonionic interpretation of the two other pairs.

\subsection{Three-forms in eight variables}
For completeness let us briefly discuss the case of the degree sixteen $SL(V_8)$-invariant of $\wedge^3V_8$, the 
last prehomogeneous space of skew-symmetric three-forms. This invariant is given by the equation 
of the projective dual of the Grassmannian $G(3,V_8)$.   This case is associated
with the grading of length seven of $\fe_8$ defined by the simple root $\alpha_2$:
$$\fe_8 = V_8^\vee \oplus \wedge^2V_8\oplus \wedge^3V_8^\vee \oplus \fgl(V_8)\oplus \wedge^3V_8\oplus \wedge^2V_8^\vee\oplus V_8.$$
(For simplicity we wrote this grading as a decomposition into $\fsl(V_8)$-modules.)
The adjoint action induces maps $Sym^{2k}\fg_1\ra Hom(\fg_{-k},\fg_k)$, which for $k=2$ and $k=3$
yield the following:
$$ P : Sym^4(\wedge^3V_8) \lra Hom (\wedge^2V_8, \wedge^2V_8^\vee), $$
$$ Q : Sym^6(\wedge^3V_8) \lra Hom (V_8^\vee, V_8). $$
Putting together $P$ and $Q$, and using the contraction map $\wedge^2V_8^\vee\otimes V_8\ra V_8^\vee$ and its dual, we get a morphism
$$ R : Sym^{10}(\wedge^3V_8) \lra Hom (V_8, V_8^\vee). $$
The computations made by Kimura in \cite[Example 2.7]{kimurabook} imply the following result:

\begin{theo} \label{16mf}
The pair of symmetric  morphisms $(Q, R)$ is a matrix factorization
of the degree sixteen semi-invariant $J_{16}$ of $\wedge^3 V_8$.
\end{theo}

\medskip

\bibliographystyle{alpha}

\begin{thebibliography}{}

\end{thebibliography}


\begin{thebibliography}{Aa}

\bibitem{abuaf} Abuaf R., {\it On quartic double fivefolds and the matrix 
factorizations of exceptional quaternionic representations}, arXiv:1709.05217.

\bibitem{baez} Baez J., {\it The octonions},  Bull. Amer. Math. Soc. {\bf 39} (2002), 145--205.

\bibitem{beauville} Beauville A., {\it Determinantal hypersurfaces},
Michigan Math. J. {\bf 48} (2000), 39--64. 

\bibitem{bertin} Bertin J., {\it 
Clifford algebras and matrix factorizations},
Adv. Appl. Clifford Algebr. {\bf 18} (2008), 417--430. 

\bibitem{beh}
Buchweitz R.O.,  Eisenbud D., Herzog J., {\it 
Cohen-Macaulay modules on quadrics}, in Singularities, representation of algebras, and vector bundles (Lambrecht, 1985), 58--116, 
Lecture Notes in Math. {\bf 1273}, Springer 1987.

\bibitem{cdp} Candelas P., Derrick E., Parkes L., {\it 
Generalized  Calabi-Yau  Manifolds  and  the Mirror of a Rigid Manifold}, 
Nucl.Phys. B {\bf 407}  (1993), 115--154.

\bibitem{chevalley} Chevalley C., {\it  The algebraic theory of spinors and Clifford algebras}, in Collected works, vol. 2,
Springer 1997.

\bibitem{clerc} Clerc J.-L., {\it Special prehomogeneous vector spaces associated 
to F4, E6, E7, E8 and simple Jordan algebras of rank 3}, 
J. Algebra  {\bf 264} (2003), 98--128. 

\bibitem{dg} Deligne P.,  Gross, B., {\it
On the exceptional series, and its descendants},
C. R. Math. Acad. Sci. Paris {\bf 335} (2002),  877--881.

\bibitem{eisenbud}
Eisenbud D., {\it  Homological algebra on a complete intersection, with an application to group 
representations}, Trans. Amer. Math. Soc. {\bf 260} (1980), 35--64.

\bibitem{gyoja} Gyoja A., {\it Construction of invariants}, 
Tsukuba J. Math. {\bf 14} (1990) 437--457. 

\bibitem{kac} Kac V., {\it 
Some remarks on nilpotent orbits}, J. Algebra {\bf 64} (1980), 190–213. 

\bibitem{igusa} Igusa  J. {\it 
A classification of spinors up to dimension twelve},
Amer. J. Math. {\bf 92} (1970),  997--1028. 

\bibitem{im} Iliev A., Manivel L., {\it Fano manifolds of Calabi-Yau Hodge type}, 
J. Pure Appl. Algebra {\bf 219} (2015), 2225--2244. 

\bibitem{kv} Gatti V., Viniberghi E., {\it 
Spinors of 13-dimensional space}, 
Adv. in Math. {\bf 30} (1978),  137--155. 

\bibitem{hitchin} Hitchin N., {\it $SL(2)$ over the octonions}, arXiv:1805.02224.

\bibitem{kimura} Kimura T., {\it 
Remark on some combinatorial construction of relative invariants},
Tsukuba J. Math. {\bf 5} (1981), 101--115. 

\bibitem{kimurabook} Kimura T., 
Introduction to prehomogeneous vector spaces, 
Translations of Mathematical Monographs {\bf 215}, AMS 2003. 

\bibitem{sk} Kimura T., Sato M., {\it 
A classification of irreducible prehomogeneous vector spaces and their relative invariants},
Nagoya Math. J. {\bf 65} (1977), 1--155. 

\bibitem{kwE7} Kraskiewicz W., Weyman J., {\it Geometry of orbit closures for 
the representations associated to gradings of Lie algebras of type $E_7$}, 
arXiv:1301.0720.

\bibitem{kwE8} Kraskiewicz W., Weyman J., {\it Geometry of orbit closures for 
the representations associated to gradings of Lie algebras of type $E_8$}, preprint.

\bibitem{LMfreud} Landsberg J.M., Manivel L., {\it 
The projective geometry of Freudenthal's magic square},
J. Algebra  {\bf 239} (2001),  477--512. 
 
\bibitem{Lie} LiE,  A computer algebra package for Lie group computations, available online at 
http://wwwmathlabo.univ-poitiers.fr/$\sim$maavl/LiE/

\bibitem{cayley} Manivel L., {\it The Cayley Grassmannian},  J. Algebra {\bf 503} (2018), 277--298. 

\bibitem{doublespinor} Manivel L., {\it Double spinor Calabi-Yau varieties}, arXiv:1709.07736.

\bibitem{ulrich} Manivel L., {\it Ulrich and aCM bundles from invariant theory}, arXiv:1803.07857.

\bibitem{popov} Popov V.L., {\it 
Classification of the spinors of dimension fourteen},
Uspehi Mat. Nauk {\bf 32} (1977), 199--200. 
 
\bibitem{rosenfeld} Rosenfeld B., Geometry of Lie groups, Mathematics and its Applications {\bf 393}, Kluwer 1997.

\bibitem{sch} Schimmrigk R., {\it 
Mirror Symmetry and String Vacua from a Special Class of Fano
Varieties}, Int. J. Mod. Phys. A {\bf 11} (1996), 3049--3096.

\bibitem{st} Seidel P., Thomas R., {\it 
Braid group actions on derived categories of coherent sheaves}, 
Duke Math. J. {\bf 108} (2001),  37--108. 

\bibitem{vinberg} Vinberg E., {\it 
Classification of homogeneous nilpotent elements of a semisimple graded Lie algebra},
Trudy Sem. Vektor. Tenzor. Anal. {\bf 19} (1979), 155--177.

\end{thebibliography}

\vspace{4mm}

\noindent Roland {\sc Abuaf}

\noindent AIM 

\noindent {\tt rabuaf@gmail.com}

\medskip 

\noindent Laurent {\sc Manivel}

\noindent Institut de Math\'ematiques de Toulouse, UMR 5219

\noindent Universit\'e de Toulouse, CNRS

\noindent 
UPS, F-31062 Toulouse Cedex 9, France

\noindent {\tt manivel@math.cnrs.fr}
\end{document}